\documentclass[12pt]{article} \oddsidemargin=0in \parskip=.7ex
\begin{document}
\newcommand{\nc}{\newcommand}  \nc{\ov}{\over} \nc{\iy}{\infty}
\nc{\inv}{^{-1}}  \nc{\be}{\begin{equation}} \nc{\A}{{\rm Ai}}
\nc{\ee}{\end{equation}} \nc{\ch}{\raisebox{.4ex}{$\chi$}} 
\renewcommand{\sp}{\vskip2ex} \nc{\noi}{\noindent} \nc{\tn}{\otimes}
\nc{\twotwo}[4]{\left(\begin{array}{cc}#1&#2\\&\\#3&#4\end{array}\right)}
\nc\tr{{\rm tr}\,}

\begin{center} {\large\bf  On Asymptotics for the Airy Process}\end{center}

\begin{center}{{\bf Harold Widom}\footnote
{Supported by National Science Foundation grant DMS-0243982.}\\
{\it Department of Mathematics\\
University of California, Santa Cruz, CA 95064\\
email address: widom@math.ucsc.edu}}\end{center}

\begin{abstract}
The Airy process $t\to A(t)$, introduced by
Pr\"ahofer and Spohn, is the limiting stationary process for a polynuclear growth 
model. Adler and van Moerbeke found a PDE in the variables $s_1,\ s_2$ and $t$
for the probability ${\rm Pr}\,(A(0)\le s_1,\, A(t)\le s_2)$.
Using this they were able, assuming the truth of a certain conjecture 
and appropriate uniformity, to obtain the first few terms of an asymptotic expansion for 
this probability as $t\to\iy$, with fixed $s_1$ and $s_2$. We shall show that
the expansion can be obtained by using the Fredholm determinant representation for
the probability. The main ingredients are formulas obtained by the author and C. A. Tracy in 
the derivation of the Painlev\'e II representation for the distribution function $F_2$ plus a few others 
obtained in the same way.
\end{abstract}

The Airy process $t\to A(t)$, introduced by
Pr\"ahofer and Spohn \cite{PS} (see also \cite{J}), is the limiting stationary process for a polynuclear growth 
model. It is also the limiting process for the largest eigenvalue of a Hermitian matrix whose
entries undergo a Dyson Brownian motion. For any fixed $t$ the probability
Pr$\,(A(t)\le s)$ equals the distribution function $F_2(s)$ which was shown in \cite{tw1} to be
representable in terms of a certain Painlev\'e~II function.

In \cite{AvM} Adler and van Moerbeke found a PDE in the variables $s_1,\ s_2$ and $t$
for the probability
\[{\rm Pr}\,(A(0)\le s_1,\, A(t)\le s_2).\]
Using this they were able, assuming the truth of a certain conjecture
and appropriate uniformity, to obtain the first few terms of an asymptotic expansion for 
this probability as $t\to\iy$, with fixed $s_1$ and $s_2$. It had the form
\[c_2(s_1,s_2)\,t^{-2}+c_4(s_1,s_2)\,t^{-4}+O(t^{-6}),\]
with explicitly computed coefficients $c_2$ and $c_4$.
The form of the coefficients suggests that
the expansion might be obtained by using the Fredholm determinant representation for
the probability, and we shall show here that this is so. After a straightforward 
computation the main ingredients will be some of the formulas obtained in 
\cite{tw1} to derive the Painlev\'e II representatioin of $F_2$ plus a few others obtained in much the
same way.

We write $\ch_i$ for the function $\ch(x-s_i)$. The probability in question is the determinant of 
$I$ minus the operator
with kernel
\[\twotwo{\ch_1(x)\,{\displaystyle\int_0^\iy} \A(x+z)\,\A(y+z)\,dz\,\ch_1(y)}
{\ch_1(x)\,{\displaystyle\int_0^\iy} e^{-zt}\,\A(x+z)\,\A(y+z)\,dz\,\ch_2(y)}
{-\ch_2(x)\,{\displaystyle\int_{-\iy}^0} e^{zt}\,\A(x+z)\,\A(y+z)\,dz\,\ch_1(y)}
{\ch_2(x)\,{\displaystyle\int_0^\iy} \A(x+z)\,\A(y+z)\,dz\,\ch_2(y)}.\]
If we set
\[K(x,y)=\int_0^\iy \A(x+z)\,\A(y+z)\,dz,\]
then the above equals
\[\twotwo{\ch_1(x)\,K(x,y)\,\ch_1(y)}{0}{{\hspace{-4ex}}0}{{\hspace{-4ex}}\ch_2(x)\,K(x,y)\,\ch_2(y)}\]
\sp
\[+\twotwo{0}{{\hspace{-8ex}}\ch_1(x)\,{\displaystyle\int_0^\iy}e^{-zt} \A(x+z)\,\A(y+z)\,dz\,\ch_2(y)}
{-\ch_2(x)\,{\displaystyle\int_{-\iy}^0} e^{zt}\,\A(x+z)\,\A(y+z)\,dz\,\ch_1(y)}{{\hspace{-8ex}}0}.\]

The second operator has an asymptotic expansion
as $t\to\iy$ which, by easy estimates, is valid in trace norm. The upper-right corner is 
\[t^{-1}\,\ch_1\,\A\tn \A\,\ch_2+t^{-2}(\ch_1\,\A'\tn \A\,\ch_2+\ch_1\,\A\tn \A'\,\ch_2)\]
\[+t^{-3}(\ch_1\,\A''\tn \A\,\ch_2+2\,\ch_1\,\A'\tn \A'\,\ch_2+\ch_1\,\A\tn \A''\,\ch_2)+\cdots,\]
while the lower-left corner is
\[-t^{-1}\,\ch_2\,\A\tn \A\,\ch_1+t^{-2}(\ch_2\,\A'\tn \A\,\ch_1+\ch_2\,\A\tn \A'\,\ch_1)\]
\[-t^{-3}(\ch_2\,\A''\tn \A\,\ch_1+2\,\ch_2\,\A'\tn \A'\,\ch_1+\ch_2\,\A\tn \A''\,\ch_1)+\cdots,\]
Here $f\tn g$ denotes either the function $f(x)\,g(y)$ or the operator with this kernel.

We will factor out $I$ minus the main operator
\[\twotwo{\ch_1(x)\,K(x,y\,\ch_1(y)}{0}{{\hspace{-4ex}}0}
{{\hspace{-4ex}}\ch_2(x)\,K(x,y\,\ch_2(y)},\]
which has determinant $F_2(s_1)\,F_2(s_2)$. For the resulting operator we use the notations
from \cite{tw1,tw2},
\[(I-\ch\,K\,\ch)\inv\,\ch\,\A=Q,\ \ (I-\ch\,K\,\ch)\inv\,\ch\,\A'=P,
\ \ (I-\ch\,K\,\ch)\inv\,\ch\,\A''=Q_1.\]
Here $\ch$ could be $\ch_1$ or $\ch_2$ and we use the notations $Q(s_1),\ Q(s_2)$, etc.
to distinguish them.\footnote{In \cite{tw1,tw2} we defined $Q(s)=(I-K\,\ch)\inv\A$,
etc., but the functions agree on $(s,\,\iy)$ so in the end that will not matter.}

After factoring out on the left $I$ minus the main operator we are left with $I$
minus the operator with matrix kernel $T$ whose diagonal entries are zero, whose upper-right
corner is 
\pagebreak
\[T_{12}:=t\inv\,Q(s_1)\tn \A\,\ch_2+t^{-2}(P(s_1)\tn \A\,\ch_2+Q(s_1)\tn \A'\,\ch_2)\]
\[+t^{-3}(Q_1(s_1)\tn \A\,\ch_2+2\,P(s_1)\tn \A'\,\ch_2+Q(s_1)\tn \A''\,\ch_2)
+\cdots,\]
and whose lower-left corner is
\[T_{21}:=-t\inv\,Q(s_2)\tn \A\,\ch_1+t^{-2}(P(s_2)\tn \A\,\ch_1+Q(s_2)\tn \A'\,\ch_1)\]
\[-t^{-3}(Q_1(s_2)\tn \A\,\ch_1+2\,P(s_2)\tn \A'\,\ch_1+Q(s_2)\tn \A''\,\ch_1)
+\cdots.\]

We have
\[\det\,(I-T)=\det\,(I-T_{12}T_{21})=e^{\tr\log (I-T_{12}T_{21})}\]
\[=-\tr\, T_{12}T_{21}+{1\ov 2}\Big((\tr\,T_{12}T_{21})^2-\tr\,(T_{12}T_{21})^2\Big)+\cdots.\]

To evaluate traces of products we use the fact $(f\tn g)\,(h\tn k)=(g,h)\,f\tn k$, whose 
trace is $(g,h)\,(f,k)$, and from \cite{tw1,tw2} the notations
\[u=(Q,\,A\,\ch),\ v=(P,\, \A\,\ch)=(Q,\,\A'\,\ch),\ w=(P,\, \A'\,\ch),
\ u_1=(Q_1,\,\A\,\ch)=(Q,\,\A''\,\ch).\]
We find that 
\[\tr\, T_{12}T_{21}=-u(s_1)\,u(s_2)\,t^{-2}
-[v(s_1)v(s_2)+2u_1(s_1)u(s_2)
-w(s_1)u(s_2)+{\rm reversed}]\,t^{-4}+\cdots,\]
\[\tr\,(T_{12}T_{21})^2=u(s_1)^2u(s_2)^2\,t^{-4}+\cdots,\]
where ``reversed'' denotes the same terms as before but with $s_1$ and $s_2$ interchanged.

It follows first that $c_2(s_1,s_2)=F_2(s_1)F_2(s_2)u(s_1)\,u(s_2)$, and since
$u=F_2'/F_2$ (see below) we have $c_2(s_1,s_2)=F_2'(s_1)\,F_2'(s_2)$ in agreement with 
\cite{AvM}. We see also that $c_4(s_1,s_2)$ equals $F_2(s_1)\,F_2(s_2)$ times
\[v(s_1)\,v(s_2)+(2\,u_1(s_1)-w(s_1))\,u(s_2)
+{\rm reversed}.\]

At this point we use formulas from \cite{tw1,tw2}, which give 
representations for $u,\, v,\, w$ and $u_1$ in terms of the function $q$ which satisfies 
the P$_{II}$ equation $q''=sq+2q^3$ and
$q(s)\sim\A(s)$ as $s\to\iy$. (About half of these
already appear in \cite{tw1}.)

In the notation of formulas (2.15)--(2.18) of \cite{tw2}, we have\footnote
{The definitions of the terms appearing on the right sides are $q=Q(s+),\ p=P(s+),\ q_1=Q_1(s+)$.}
\[u'=-q^2,\ v'=-pq,\ w'=-p^2,\  u_1'=-q_1\,q,\]
and by formula (2.12) of \cite{tw2} with $x$ replaced by $s$
\[q_1(s)=sq-vq+up.\]

By the formula for $u'$, since $u(+\iy)=0$,
\[u=\int_s^\iy q(x)^2.\]
(Thus $F_2'/F_2=u$, as stated above.) Eventually everything will be expressed in terms of 
$q$ and $u$.

We have $v'=-pq$ and formula (3.1) of \cite{tw1}, which says $q'=p-qu$, gives
\[v'=-(q'+qu)\,q=-{1\ov2}(q^2)'+{1\ov2}(u^2)',\]
so
\[v=-{1\ov2}q^2+{1\ov2}u^2.\] 

Next,
\[u_1'=-q_1q=-(s\,q-vq+up)\,q=-sq^2+vq^2-qu(q'+qu)\]
and
\[-w'=(q'+qu)^2,\]
which give
\[(2\,u_1-w)'=-2sq^2+2vq^2+(q')^2-q^2u^2.\]
By the formula for $v$ this is
\[-2sq^2-q^4+(q')^2,\]
and so
\[2\,u_1-w=\int_s^\iy (2x\,q(x)^2+q(x)^4-q'(x)^2)\,dx.\]

Putting these together we find that $c_4(s_1,s_2)$ equals
$F_2(s_1)\,F_2(s_2)$ times
\[{1\ov4}u(s_1)^2\,u(s_2)^2+{1\ov4}q(s_1)^2\,q(s_2)^2-{1\ov2}q(s_1)^2\,u(s_2)^2\]
\[+\int_{s_1}^\iy (2x\,q(x)^2+q(x)^4-q'(x)^2)\,dx\,
\cdot\, u(s_2)+{\rm reversed}.\]

This looks a little different from the formula of \cite{AvM}, because there the
last integral is (we change $s_1$ back to $s$)
\[\int_s^\iy (2(s-x)\,q(x)^2-q(x)^4+q'(x)^2)\,dx.\]
But if we take $d^2/ds^2$ of the two integrals we see, using
the equation $q''=sq+2q^3$ satisfied by $q$, that the results are the
same. Since the integrals and their derivatives both vanish at $s=+\iy$ the integrals
are equal.

\end{document}